\documentclass[12pt]{article}

\usepackage[centertags]{amsmath}
\usepackage{amsfonts}
\usepackage{amsthm}
\usepackage{newlfont}
\usepackage{amscd}
\usepackage{amsgen}
\usepackage{amssymb}
\usepackage{stmaryrd}
\usepackage{amssymb,amsmath}
\usepackage{amscd}
\usepackage{enumerate}

\newlength{\defbaselineskip} 
\theoremstyle{plain}
\newtheorem{thm}{Theorem}[section]
\newtheorem{cor}[thm]{Corollary}
\newtheorem{con}[thm]{Conjecture}
\newtheorem{df}[thm]{Definition}
\newtheorem{lema}[thm]{Lemma}
\newtheorem{obs}[thm]{Proposition}
\newtheorem{exm}[thm]{Example}

\newtheorem{rem}[thm]{Remark}

\theoremstyle{definition} 
\theoremstyle{definition} \newtheorem{ex}{Example}[section] %

 \numberwithin{equation}{section}

\def\z{\mathbb{Z}}

\def\o{\mathcal{O}}

\def\ob{\begin{obs}}
\def\kob{\end{obs}}
\def\dow{\begin{proof}}
\def\kdow{\end{proof}}
\def\kwadrat{\hfill$\square$}
\def\tw{\begin{thm}}
\def\ktw{\end{thm}}
\def\hip{\begin{con}}
\def\khip{\end{con}}
\def\lem{\begin{lema}}
\def\klem{\end{lema}}
\def\ex{\begin{exm}}
\def\prog{\begin{pr}}
\def\kprog{\end{pr}}
\def\wn{\begin{cor}}
\def\kwn{\end{cor}}
\def\uwa{\begin{rem}}
\def\kuwa{\end{rem}}
\def\kex{\end{exm}}
\def\dfi{\begin{df}}
\def\kdfi{\end{df}}

\begin{document}

\title{Derived category of toric varieties with Picard number three}

\author{Arijit Dey \\
Micha\l\hspace{2pt}  Laso\'{n} \\
Mateusz Micha\l ek}


\maketitle
\begin{abstract}
We construct a full, strongly exceptional collection of line
bundles on the variety $X$ that is the blow up of the
projectivization of the vector bundle $\o_{\mathbb P^{n-1}}\oplus
\o_{\mathbb P^{n-1}}(b_1)$ along a linear space of dimension
$n-2$, where $b_1$ is a non-negative integer.
\end{abstract}

\section{Introduction}
Let $X$ be a smooth projective variety over the field of complex numbers
$\mathbb C$ and let $D^b(X)$ be the derived category of bounded complexes of 
coherent sheaves of $\mathcal O_X$-modules. This category is an important 
algebraic invariant of $X$. In order to understand the derived category 
$D^b(X)$ one is interested in knowing a {\it strongly exceptional} finite 
collection of objects in $D^b(X)$ that generates the derived category $D^b(X)$. 
 
The notion of ''strongly exceptional'' collection was first 
introduced by Gorodentsev and Rudakov \cite{gr} in order to study vector 
bundles on $\mathbb P^n$. An {\it exceptional} 
collection $\{F_0,\,F_1,\,\cdots,\,F_{m}\}$ of sheaves gives a functor  
$F_{E}$ from the category of coherent sheaves $Coh(X)$ to 
the derived category $D^b(\mathcal A-module)$ of $\mathcal A$-modules, where $E\,=\,\oplus_{i=0}^{m} F_i$ and $\mathcal A\,=\,Hom(E,E)$.
The functor $F_E$ is extendable to the derived functor $D^b(F_E)$ from 
$D^b(X)$ to $D^b(\mathcal A-module)$. In \cite{bondal1} Bondal proved that if $\{F_0,\,F_1,\,\cdots,\,F_{m}\}$ is a {\it full strongly exceptional} collection then the functor 
$D^b(F_E)$ is an equivalence of categories. 
The existence of a {\it full strongly exceptional collection} $\{F_0,\,F_1,\,\cdots,\,F_{m}\}$ of 
coherent sheaves on a smooth projective variety puts a strong condition on $X$ 
namely the Grothendieck group 
$K^{0}(X)$ is isomorphic to $\mathbb Z^{m+1}$. In general for an arbitrary 
variety $X$, $K^0(X)$ has torsion 
but $K^0$ of a toric variety is a finitely generated free abelian group. So 
it is interesting to search for {\it full strongly exceptional collections} of 
sheaves in case of toric variety. 
For a smooth complete toric variety Kawamata \cite{kawamata} proved that the derived category $D^b(X)$ has a full collection of exceptional objects. In his 
collection, the objects are sheaves rather than line bundles and the collection is only exceptional and not strongly exceptional. 
For a smooth complete toric variety $X$, there is a well known construction 
due to Bondal which gives a (finite) full collection of line bundles of $D^b(X)$. In general 
Bondal's collection of line bundles need not be a strongly exceptional 
collection but one "hopes" that for huge families of toric varieties we will be able to choose a subset and order it in such a 
way that it becomes a full strongly exceptional collection. 
In \cite{king} King made the following conjecture:

\hip For any smooth, complete toric variety $X$ there
exists a full, strongly, exceptional collection of line bundles.
\footnote{Originally this conjecture was made in terms of
existence of titling bundles whose direct summands are line
bundles but it is easy to see that they are equivalent, see
\cite{rosa-costa-tilting}} \khip This conjecture was disproved by
Hille and Perling, in \cite{Hille-Perling} they gave an example of
a smooth, complete toric surface that does not have a full
strongly exceptional collection of line bundles. The conjecture
was reformulated by Costa and Mir\'o-Roig:

\hip For every smooth, complete Fano toric variety there exists a
full strongly exceptional collection of line bundles. \khip This
conjecture is still open and is supported by many numerical
evidences. It has an affirmative answer when the Picard number of $X$
is less than or equal to two. When the Picard number is one it is
a full strongly exceptional collection easy to see that $X$ is isomorphic to projective space $\mathbb
P^r$ and the collection $(\mathcal O,\,\mathcal
O(1),\,\cdots,\,\mathcal O(r))$ is a full strongly exceptional collection (this is a consequence of Beilinson's theorem). When the Picard number of $X$
is two, the above question has affirmative answer and this is due
to Costa and Mir\'o-Roig \cite{rosa-costa-tilting}. When the Picard 
number is $3$ the question is not fully resolved.    

Toric varieties with Picard number three are completely classified
by Batyrev \cite{baty} in terms of its {\it primitive
collections}, he showed that any toric variety with Picard number
$3$ has $3$ or $5$ primitive collections. Toric varieties with $3$
primitive collections are isomorphic to a projectivization of a
decomposable bundle over a smooth toric variety $W$ of a smaller
dimension with Picard number $2$, hence by \cite{rosa-costa-tilting} 
we have an affirmative answer to the conjecture. When the number of 
primitive collections is $5$ the conjecture is still open. There are 
some partial results known in this case, for example recently
Costa and Mir\'o-Roig \cite{lcmr} answered the above conjecture
affirmatively when $X$ is a blow up of $\mathbb P^{n-r} \times
\mathbb P^r$ along a multilinear subvariety of codimension $2$.
Motivated by this result, R. M. Mir\'o-Roig and L. Costa (in the
meeting P.R.A.G.MAT.I.C' 09) suggested us to investigate this
question for a large family of toric varieties parameterized by
positive integers $b_1$, $n$. In this note we consider a toric
variety $X$ which is a blow up of the projectivization of the vector
bundle $\o_{\mathbb P^{n-1}}\oplus \o_{\mathbb P^{n-1}}(b_1)$ on
$\mathbb P^{n-1}$ along a linear space of dimension $n-2$, 
where $b_1$ is a positive integer. We are able to answer the
conjecture affirmatively for this family of toric varieties (see Theorem \ref{glowne}).
Note that not all our varieties are Fano, in fact $X$ is Fano 
if $b_1 < n-1$.

We outline the structure of this paper. In \S 2, we briefly review
the notions of strongly exceptional collection of sheaves and few
basic facts about toric varieties which will be needed later on.
In \S 3, we recall Batyrev's classification of toric varieties and
we describe the family of toric varieties which we are interested in, in
terms of fans and its primitive relations. In \S 4 we determine
explicitly a full strongly exceptional collection of line bundles
for this family of toric varieties.

\section*{Acknowledgements}
We would like to thank very much L. Costa and R.M. Mir\'o-Roig for
introducing us to the subject. We are also
grateful to them for sharing their many useful ideas and meticulous 
reading of a preliminary draft of this paper.
\section{Preliminaries}
The goal of this section is to fix the notation and basic facts
that we will use through this paper. We start by recalling the
notions of exceptional sheaves, exceptional collection of sheaves,
strongly exceptional collection of sheaves and full strongly
exceptional collection of sheaves. Let $X$ be a smooth projective variety 
over $\mathbb C$.
\begin{df} \hskip 0pt
\begin{enumerate}
\item{} A coherent sheaf $F$ on $X$ is {\bf exceptional} if
$\text{Hom}(F,F)= \mathbb C$ and \text{Ext}$\,^i_{\o_X}(F,F)
\,=\,0$ for $i > 0$. \item{} An ordered collection
$(F_0,\,F_1,\,\cdots,\,F_{m})$ of coherent sheaves on $X$ is an
{\bf exceptional collection} if each sheaf $F_i$ is exceptional
and \\$\text{Ext}^i_{\o_X}(F_k,F_j)=0$ for $j<k$ and $i\geq 0$.
\item{} An exceptional collection $(F_0,\,F_1,\,\cdots,\,F_{m})$
of coherent sheaves on $X$ is a {\bf strongly exceptional
collection} if $\text{Ext}^i_{\o_X}(F_j,F_k)=0$ for $j \leq k$ and
$i \geq 1$. \item{} A (strongly) exceptional collection
$(F_0,\,F_1,\,\cdots,\,F_{m})$ of coherent sheaves on $X$ is a
{\bf full (strongly) exceptional collection } if it generates the
bounded derived category $D^b(X)$ of $X$ i.e. the smallest
triangulated category containing $\{F_0,\,F_1,\,\cdots,\,F_m\}$ is
equivalent to $D^b(X)$.
\end{enumerate}
\end{df}

\section{Toric varieties with Picard number three}
In this section we introduce notation and facts concerning toric varieties that we
use in our paper. An $n$ dimensional toric variety $X$ is a smooth,
projective variety containing an $n$ dimensional torus $T$ ($n$
copies of $\mathbb C^*$) together with an action on $X$ and
characterized by a fan $\Sigma$ of strongly convex polyhedral
cones in $N \otimes_{\mathbb Z} \mathbb R$, where $N$ is a lattice
$\mathbb Z^n$. We denote the $\mathbb Z$-basis of $N$ by
$e_1\,,\cdots\,,e_n$ and by $e_1^*\,,\cdots, \,e_n^*$ its dual basis 
in $M:=\,Hom_{\mathbb Z}(N,\mathbb Z)$. For every one dimensional cone
$\sigma\,\in\,\Sigma$ there is a unique generator $v \in N$ such
that $\sigma \cap N\,=\,\mathbb Z_{\ge0}\cdot v$, it is called the ray
generator. The set of all ray generators is denoted by $R$. To
each ray generator $r \in R$ one can associate a toric divisor
$D_r$ (see also \cite{Fulton}). If the number of toric divisors is
$m$ then the Picard number of $X$ is $m-n$ where $n$ is the dimension
of $X$. The anticanonical divisor $-K_X$ is given by
$-K_X\,=\,\sum_{r\in R}D_r$. We say that $X$ is Fano if $-K_X$ is ample. 

Smooth, complete
toric varieties with Picard number three have been classified by
Batyrev in \cite{baty} according to their primitive relations. Let
$\Sigma$ be a fan in $N=\z^n$.

\dfi We say that a subset $P\subset R$ is a primitive collection
if it is a minimal (with respect to inclusion) subset of $R$ which
does not span a cone in $\Sigma$. \kdfi

In other words a primitive collection is a subset of ray
generators, such that all together they do not span a cone in
$\Sigma$ but if we remove any generator, then the rest spans a
cone that belongs to $\Sigma$. To each primitive collection
$P=\{x_1,\dots, x_k\}$ we associate a primitive relation. Let
$w=\sum_{i=1}^k x_i$ and $\sigma \in \Sigma$ be the cone of the
smallest dimension that contains $w$. Let $y_1,\dots,y_s$ be
the ray generators of this cone. The toric variety of $\Sigma$ was
assumed to be smooth, so there are unique positive integers
$n_1,\dots, n_s$ such that
$$w=\sum_{i=1}^s n_iy_i.$$
\dfi For each primitive collection $P=\{x_1,\dots,x_k\}$ the linear 
relation:
$$x_1+\dots+x_k-n_1y_1-\dots-n_sy_s=0$$
is called the primitive relation (associated to $P$). \kdfi Using
the results of \cite{grun} and \cite{odap} Batyrev proved in
\cite{baty} that for any smooth, complete $n$ dimensional fan with $n+3$
generators its set of ray generators can be partitioned into $l$ non-empty
sets $X_0,\dots,X_{l-1}$ in such a way that the primitive
collections are exactly sums of $p+1$ consecutive sets $X_i$ (we
use a circular numeration, that is we assume that $i\in\z/l\z$), where
$l=2p+3$. Moreover $l$ is equal to $3$ or $5$. The number $l$ is
of course the number of primitive collections. In the case $l=3$
the fan $\Sigma$ is a splitting fan (that is any two primitive
collections are disjoint). These varieties are well characterized,
and we know much about full strongly exceptional collections of
line bundles on them. The case of five primitive
collections is much more complicated and is our object of study. For
$l=5$ we have the following result of Batyrev \cite[Theorem 6.6]{baty}.

\tw Let $Y_i=X_i\cup X_{i+1}$, where $i\in\z/5\z$,
$$X_0=\{v_1,\dots,v_{p_0}\},\quad X_1=\{y_1,\dots,y_{p_1}\},\quad X_2=\{z_1,\dots,z_{p_2}\},$$
$$X_3=\{t_1,\dots,t_{p_3}\},\quad X_4=\{u_1,\dots,u_{p_4}\},$$
where $p_0+p_1+p_2+p_3+p_4=n+3$. Then any $n$-dimensional fan
$\Sigma$ with the set of generators $\bigcup X_i$ and five
primitive collections $Y_i$ can be described up to a symmetry of
the pentagon by the following primitive relations with nonnegative
integral coefficients $c_2,\dots,c_{p_2},b_1,\dots,b_{p_3}$:
$$v_1+\dots+v_{p_0}+y_1+\dots+y_{p_1}-c_2z_2-\dots-c_{p_2}z_{p_2}-(b_1+1)t_1-\dots-(b_{p_3}+1)t_{p_3}=0,$$
$$y_1+\dots+y_{p_1}+z_1+\dots+z_{p_2}-u_1-\dots-u_{p_4}=0,$$
$$z_1+\dots+z_{p_2}+t_1+\dots+t_{p_3}=0,$$
$$t_1+\dots+t_{p_3}+u_1+\dots+u_{p_4}-y_1-\dots-y_{p_1}=0,$$
$$u_1+\dots+u_{p_4}+v_1+\dots+v_{p_0}-c_2z_2-\dots-c_{p_2}z_{p_2}-b_1t_1-\dots-b_{p_3}t_{p_3}=0.$$
\ktw

In our case we will be interested in varieties $X$ with Picard number
three that have the following sets $X_i$:
\begin{eqnarray}\label{pc} 
X_0=\{v_1,\dots,v_{n-1}\},\, X_1=\{y\},\, X_2=\{z\},\, X_3=\{t\},\, X_4=\{u\}
\end{eqnarray} 
So, from now on let us denote by $X$ a smooth toric variety with Picard number 3 and primitive collections $X_0\cup X_1$, $X_1\cup X_2,\dots,X_4\cup X_0$.
We see that the cone $<v_1,\dots,v_{n-1},z>$ is in our fan, so as
the variety is smooth, these ray generators form a basis
$(e_1,\dots,e_n)$ of a lattice. In this basis, using the above primitive
relations we see that all the considered ray generators are of
the following form:
\begin{eqnarray} \label{pr}
\begin{cases}
v_1=e_1,\dots,v_{n-1}=e_{n-1},\\
t=-e_n,\\
z=e_n, \\
u=-e_1-\dots-e_{n-1}-be_n, \\
y=-e_1-\dots-e_{n-1}-(b+1)e_n.
\end{cases}
\end{eqnarray}
One can see that for any fixed dimension $n$ we obtain an infinite
number of smooth toric varieties parameterized by $b=b_1\geq 0$, but
only a finite number of them is Fano, namely for $b< n-1$ (because the sum of coefficients in each primitive relation has to be positive).

We need following basic facts about divisors
on toric varieties. To each ray generator $r \in R$ we can  
associate a divisor $D_r$ \cite{Fulton}. The relations (in the
Picard group) among the divisors are given by the following equations:
$$\sum_{r\in R} e_i^*(r)D_r=0,$$
what is
$$D_{v_1}-D_u-D_y=0,\dots,D_{v_{n-1}}-D_u-D_y=0,$$
$$D_z-D_t-bD_u-(b+1)D_y=0.$$

\lem
The above linear relations imply that $Pic(X)\cong \z^3=<D_v,D_y,D_t>$, i.e. each divisor can be uniquely written in a form:
$$eD_v+fD_y+gD_t,$$
where $D_v$ is any fixed $D_{v_i}$ (they are all linearly equivalent).
\klem
\dow From the relations above it is obvious that each divisor is linearly equivalent to a divisor of this form. Let us assume that it has two such presentations. It means that they have to be linearly dependent:
$$eD_v+fD_y+gD_t=e'D_v+f'D_y+g'D_t+i(D_v-D_u-D_y)$$$$+j(D_z-D_t-bD_u-(b+1)D_y).$$
Since $D_z$ occurs only on the right hand side in the above equality we have $j=0$. Once $j=0$ the divisor $D_u$ occurs only on the right hand side so $i=0$ and we get uniqueness.\kdow

\section{Main theorem}

In this section we prove that for smooth toric projective
varieties $X$ with Picard number $3$ with sets of generators
$$X_0=\{v_1,\dots,v_{n-1}\},\quad X_1=\{y\},\quad X_2=\{z\},\quad X_3=\{t\},\quad X_4=\{u\},$$
in the situation described in the previous section there exists a
full strongly exceptional collection of line bundles in the
derived category.

We proceed in several steps. First by
pushing forward a trivial line bundle by a Frobenius morphism we
obtain a vector bundle that splits into the direct sum of line bundles which by Bondal's result \cite{Ober} generate $ D^b(X)$. 
We can calculate the set $C$ of these line bundles explicitly using the algorithm described in \cite{thom}. Then, we choose an ordered subset $C' \subset C$ and we prove that $C'$ is strongly exceptional. Finally using Koszul complexes we prove that $C$ and $C'$ generate the same category, hence $C'$ is also full.

\subsection{Full collection}\label{fullcoll}
We fix a prime integer $p>>0$. Let $F:X\rightarrow X$ be the $p$-th Frobenius morphism of our toric
variety $X$, that is an extension of a morphism:
\[\begin{array}{ll}
F:&T\rightarrow T,\\
&t\rightarrow t^p
\end{array}\]
where $T$ is the torus of $X$. Using the results of \cite{thom} we
can calculate the split of the push forward $F_*(\mathcal{O}_X)$.
We will use similar notation as used in \cite{thom}. Let us recall the algorithm. We fix a basis of $N$. Let 
$$R=k[(X^{e_1^*})^{\pm 1},\dots,(X^{\hat e_n^*})^{\pm 1}],$$
be the coordinate ring of the torus $T$.
To each cone $\chi_i\subset N$ of maximal dimension we associate a matrix $A_i$ whose rows are ray generators in the chosen basis. Let $B_i=A_{i}^{-1}$ and 
$C_{ij}=B_{j}^{-1}B_{i}$. Let $w_{ij}=(w_{ij}^1,\dots,w_{ij}^n)$ be the $j$-th column of the matrix $B_i$. To each maximal cone $\chi_i\subset N$ one can also associate an open affine subvariety $U_{\chi_i}$ with the coordinate ring 
$$R_i=k[X_{i1},\dots,X_{in}]\subset R,$$
where we use the notation $X_{ij}=X^{w_{ij}}=(X^{e_1^*})^{w_{ij}^1}\hdots (X^{e_n^*})^{w_{ij}^n}$.
If we consider two cones $\chi_i,\chi_j\subset N$ then $\chi_i\cap\chi_j$ is a face of $\chi_i$. Using \cite{Fulton} Proposition 2, p.~13 we see that there is a monomial $M_{ij}$ such that $(R_i)_{M_{ij}}$ is the coordinate ring of $\chi_i\cap\chi_j$. Let
$$I_{ij}=\{v=(v_1,\dots,v_n): X_i^v\text{ is a unit in }(R_i)_{M_{ij}}\}.$$
 Let us also define the set
$$P_p=\{(g_1,\dots,g_n):0\leq g_i< p)\}.$$
For $w\in I_{ji}$ one can define maps:
$$h^w_{ijp}:P_p\rightarrow I_{ji},$$
$$r^w_{ijp}:P_p\rightarrow P_p$$
determined for $g\in P_p$ by the equality
\begin{equation}\label{eqn1}
C_{ij}g+w=ph^w_{ijp}(g)+r^w_{ijp}(g).
\end{equation}
Let us fix a Cartier divisor $D=\{(U_{\chi_i},X_i^{u_i})\}$ and a line bundle $L\cong \o(D)$. From \cite[Lemma 4]{thom} one gets, for each $g\in P_p$ and each cone $\chi_l$ a $T$-Cartier divisor $D_g=\{(U_{\chi_i},X_i^{g_i})\}$, where $g_i=h^{u_{li}}_{lip}(g)$, that is independent from the choice of $l$. Moreover 
by \cite[Theorem 1]{thom} taking all $g\in P_p$ one gets line bundles that form a split of the push forward by the Frobenius morphism $F_*(L)$.
In our case the algorithm simplifies.

 Let us consider
three matrices:
\[
A_0=Id_n,A_1=
\left[
\begin{array}{cccc}
&&&0\\
&Id_{n-1}&&\vdots\\
&&&0\\
0&\hdots&0&-1\\

\end{array}
\right],
A_2=\left[
\begin{array}{ccccc}
&&&0&0\\
&Id_{n-2}&&\vdots&\vdots\\
&&&0&0\\
-1&\hdots&-1&-1&-b\\
-1&\hdots&-1&-1&-b-1\\

\end{array}
\right].
\]
The matrices above correspond to following cones:
$$\sigma_0=<v_1,\dots,v_{n-1},z>,\hskip 6pt \sigma_1=<v_1,\dots,v_{n-1},t>,\hskip 6pt \sigma_2=<v_1,\dots,v_{n-2},u,y>.$$
From \ref{eqn1} we get:
$$C_{0j}=A_jA_0^{-1}=A_j.$$
Let $g\in P_p$, as we are pushing forward trivial line bundle we want to 
calculate $h^0_{0jp}$ and $r^0_{0jp}$ that satisfy:
$$A_jg=ph^0_{0jp}(g)+r^0_{0jp}(g),$$
where $r^0_{0jm}\in\{0,\dots, p-1\}^n$. As $A_0g=p\cdot 0+g$, we see that
$D_g$ as a Cartier divisor on $X$ is given by $1$ on
$U_{\sigma_0}$. On $U_{\sigma_1}$ the divisor is given by:
\[\begin{cases}
1&\text{ if }g_n= 0\\
X^{-w_{in}}&\text{ if }g_n\neq 0,
\end{cases}
\]
and on $U_{\sigma_2}$ by:
\[\begin{cases}

X^{-sw_{j(n-1)}}X^{-sw_{jn}}&\text{ if }g_n=0 \\
X^{-sw_{j(n-1)}}X^{-sw_{jn}}\text{ or
}X^{-sw_{j(n-1)}}X^{-(s+1)w_{jn}}&\text{ if }g_n\neq 0,
\end{cases}
\]
where $s=g_1+\dots+g_{n-1}+bg_n$.
 We see that for $p>>0$, $F_*(\mathcal{O}_X)$ splits into the
direct sum of line bundles, all of which belong to one of the following three
subsets:
$$B_1=\{\o(-qD_u-(q+1)D_y-D_t):q=0,\dots,n-1+b \}$$
$$B_2=\{\o(-qD_u-qD_y-D_t):q=1,\dots,n-1+b \}$$
$$B_3=\{\o(-qD_u-qD_y):q=0,\dots,n-1\}.$$

\ob\label{full} With the above notation the line bundles from the set $B_1\cup B_2\cup B_3$
generate the derived category $D^{b}(X)$.\kob
\dow
This is a direct consequence of Bondal's result from \cite{Ober}, that the 
split of the push forward of a trivial bundle by the $m$-th Frobenius morphism generates the derived category $D^{b}(X)$ for $m$ sufficiently large.
\kdow
\subsection{Forbidden subsets}

In this subsection we want to characterize acyclic line bundles on $X$ i.e. line bundles whose higher cohomologies vanishes. We will use this
characterization to check if $Ext^i(L,M)=H^i(L^\vee\otimes M)=0$ for $i>0$.

Let $\Sigma$ be an arbitrary fan in $N=\mathbb Z^n$ with the set
of ray generators $x_1,...,x_m$ and $\mathbb P_{\Sigma}$ be the toric variety 
associated to the fan $\Sigma$. For
$I\subset\{1,\dots,m\}$ let $C_I$ be the simplicial complex
generated by sets $J\subset I$ such that $\{x_i:i\in J\}$ generates the 
cone in $\Sigma$ and for $r=(r_i:i=1,\dots,m)$. Let us define
$Supp(r):=C_{\{i:\;r_i\geq 0\}}$.

From the result of Borisov and Hua \cite{bohu} we have the
following:

\ob\label{cohogen} The cohomology $H^j(\mathbb P_{\Sigma},L)$ is
isomorphic to the direct sum over all $r=(r_i:i=1,\dots,m)$ such
that $\o(\sum_{i=1}^m r_iD_{x_i})\cong L$ of the $(n-j)$-th reduced
homology of the simplicial complex $Supp(r)$.\kob

\dfi A line bundle $L$ on $\mathbb P_{\Sigma}$ is called acyclic if 
$H^i(\mathbb P_{\Sigma},L)=0$ for $i\geq 1$.
\kdfi

\dfi A proper subset $I$ of $\{1,\dots,m\}$ is called a forbidden set if the
simplicial complex $C_I$ has nontrivial reduced homology.\kdfi

From Proposition \ref{cohogen} we have the following
characterization of acyclic line bundles

\ob A line bundle $L$ on $\mathbb P_{\Sigma}$ is acyclic if it is not isomorphic to none
of the following line bundles
$$\o(\sum_{i\in I}r_iD_{x_i}-\sum_{i\not\in I}(1+r_i)D_{x_i})$$
where $r_i\geq 0$ and $I$ is a proper forbidden subset of
$\{1,\dots,m\}$. \kob

Hence to determine which bundles on $\mathbb P_{\Sigma}$ are
acyclic it is enough to know which sets $I$ are forbidden.

In case of simplicial complex $C_I$ on the set of vertices $I$ we
also define a primitive collection as a minimal subset of vertices
that do not form a simplex. A complex $C_I$ is determined by its
primitive collections, namely it contains simplices (subsets of
$I$) that contain none of primitive collections.

In case of our variety (described at the beginning of this
section) we have $C_I=\{J\subset I: \widehat{Y_i}:=\{j:x_j\in
Y_i\}\nsubseteq J$ for $i=1,\dots,5\}$, since $Y_i$ are primitive
collections. So sets $\widehat{Y_i}$ are primitive collections in
the simplicial complex. The only difference between sets
$\widehat{Y_i}$ and $Y_i$ is that the first one is the set of
indices of rays in the second one, so in fact they could be even
identified. For our convenience we also define sets
$\widehat{X_i}:=\widehat{Y_i}\cap\widehat{Y_{i-1}}$ which are
similarly sets of indices of sets $X_i$.

\lem A primitive collection is a forbidden subset. \klem \dow Let
$I$ be a primitive collection with $k$ elements. The chain complex
of $C_I$ is as follows
$$0\rightarrow\mathbb C^{{k}\choose{k-1}}\rightarrow\mathbb C^{{k}\choose{k-2}}\rightarrow...\rightarrow\mathbb C^{{k}\choose{2}}\rightarrow\mathbb C^{{k}\choose{1}}\rightarrow\mathbb C\rightarrow 0,$$
which is not exact because the Euler characteristic is nonzero.
\kdow

\lem A sum of two consecutive primitive collections is a forbidden
subset. \klem \dow Let
$I=\widehat{Y_{i}}\cup\widehat{Y_{i+1}}=\widehat{X_i}\cup\widehat{X_{i+1}}\cup\widehat{X_{i+2}}$, $|\widehat{X_i}|=k_1$, $|\widehat{X_{i+1}}|=k_2$,
$|\widehat{X_{i+2}}|=k_3$ and $|I|=k$. Then chain complex of $C_I$
is as follows
$$0\rightarrow\mathbb C^{{{k}\choose{k-1}}-{{k_1}\choose{k-1-k_2-k_3}}-{{k_3}\choose{k-1-k_1-k_2}}}\rightarrow...\rightarrow\mathbb C^{{{k}\choose{t}}-{{k_1}\choose{t-k_2-k_3}}-{{k_3}\choose{t-k_1-k_2}}}\rightarrow...$$$$\rightarrow\mathbb C\rightarrow 0$$
which is not exact because the Euler characteristic is
nonzero.\kdow

\lem\label{notsum} If a nonempty subset $I$ is not a sum of
primitive collections, then it is not forbidden. \klem \dow
The simplicial homology of a simplicial complex is equal to the singular
homology of this complex considered as a topological space (each
simplex $D$, which is a $d$ element set, can be changed into the convex hull
of $d$ linearly independent vectors in $\mathbb R^n$ that
correspond to elements of this set). To avoid confusion with
scalars let us name elements of $\{1,\dots,m\}\supset I$ as
$\{x_1,\dots,x_m\}$. The above names are not by accident the same
as rays of a fan, because this complex as a topological space can
be realized by sum of convex hulls of sets of rays that form a
cone in $\Sigma$ and whose indices are contained in $I$.

There exists $a\in I$ such that $a$ does not belong to any
primitive collection which is contained in $I$. We can define a
homotopy
$$H:[0,1]\times C_I\rightarrow C_I$$ which for
$x=\alpha_1x_1+...+\alpha_mx_m$ ($\alpha_i\geq0$, $\sum \alpha_i=1$)  gives
$$H(t,x):=t\alpha_1x_1+...+(1-t+t\alpha_a)x_a+...+t\alpha_mx_m.$$ Of course $x\in
C_I$ means that $S_x:=\{i:\alpha_i>0\}\subset I$ and
$Y_i\nsubseteq S_x$, but then $S_x\cup\{a\}$ also
satisfies this conditions, so $H(t,x)\in C_I$ and $H$ is well
defined. It is easy to observe that $H$ is continuous
$H(0,\cdot)=x_a$ and $H(1,\cdot)=id_{C_I}$. The complex $C_I$ is homotopic to
a point, so it has trivial reduced homologies.\kdow

\lem A sum of three consecutive primitive collections is not a
forbidden subset. \klem \dow At the beginning of this proof we
should give the same remark as in the proof of lemma \ref{notsum}.
We have $I=\widehat{Y_i}\cup
\widehat{Y_{i+1}}\cup\widehat{Y_{i+2}}=\widehat{X_i}\cup\widehat{X_{i+1}}\cup\widehat{X_{i+2}}\cup\widehat{X_{i+3}}$,
so in our situation at least one of the sets $\widehat{X_{i+1}}$,
$\widehat{X_{i+2}}$ has only one element. Without loss of
generality we can assume that $\widehat{X_{i+2}}=\{x_c\}$ and also
that $\widehat{X_{i}}=\{x_{a_1},\dots,x_{a_A}\}$,
$\widehat{X_{i+1}}=\{x_{b_1},\dots,x_{b_B}\}$,
$\widehat{X_{i+3}}=\{x_{d_1},\dots,x_{d_D}\}$. Let us define the
homotopy
$$H:[0,1]\times C_I\rightarrow C_I$$ which for
$x=\alpha_{a_1}x_{a_1}+...+\alpha_{b_1}x_{b_1}+...+\alpha_{c}x_{c}+\alpha_{d_1}x_{d_1}+...+\alpha_{d_D}x_{d_D}$
gives
$$H(t,x):=x+t\alpha_{c}x_{a_1}-t\alpha_{c}x_{c}.$$
If $\alpha_{c}=0$ then $H(t,x)=x$. If $\alpha_{c}\neq0$ then
$S_{H(x,t)}\subset S_x\cup\{a_1\}$, but this set is also in our
symplicial complex $C_I$, if contrary $a_2,...,a_A,b_1,...,b_B$
are in $S_x$ so $\{b_1,...,b_B,c\}=Y_{i+1}\subset S_x$ a
contradiction. So the homotopy $H$ is well defined. It is easy to
observe that $H$ is continuous, $H(0,\cdot)=id_{C_I}$ and
$H(1,C_I)$ is a symplicial complex on vertices
$x_{a_1},...,x_{b_B},x_{d_1},...,x_{d_D}$ with only one primitive
collection $\{x_{a_1},...,x_{b_B}\}$. Hence in the same way as in
Lemma \ref{notsum} $H(1,C_I)$ can be contracted to a point
$x_{d_1}$. This shows that $C_I$ is homotopic to a point, so it has trivial
reduced homologies. \kdow

The above Lemmas match together to the following

\tw The only forbidden subsets are primitive collections, their
complements (these are exactly sums of two consecutive primitive collections) and the empty set. \ktw

This gives us that in our situation

\wn\label{acyclic} With the above notation a line bundle $L$ is acyclic if and only if it is
not isomorphic to any of the following line bundles
$$\o(\alpha_1D_{v}+\alpha_2D_{y}+\alpha_3D_{z}+\alpha_4D_{t}+\alpha_5D_{u})$$
where exactly $2,3$ or $5$ consecutive $\alpha$ are negative and
if $\alpha_1<0$ then $\alpha_1\leq-(n-1)$.\kwn \dow Since all
$D_{v_i}$ are linearly equivalent we match them together and as a
consequence $\alpha_1$ is the sum of all the coefficients of
$D_{v_i}$. \kdow

\subsection{Strongly exceptional collection}\label{strexcoll2}
We are looking for a full strongly exceptional collection. From the general theory we know that if it exists then its length should be equal to the rank of the Grothendieck group $K_0(X)$. In case of a smooth toric varieties the rank of this group is equal to the number of maximal cons in the fan. 

In our case the maximal cones are $n$ dimensional subsets of the set of all ray generators, except those subsets that contain some primitive collection. We want to calculate how many such subsets there are. First let us notice that at most $2$ elements of such subset can be contained in $X_1\cup X_2\cup X_3\cup X_4$, because otherwise it would contain a primitive collection. This means that we have got only two possibilities: 

1) There are two elements of our subset that are in this set. There are ${{n-1}\choose{n-2}}\cdot({{4\choose 2}-3})=3(n-1)$ such subsets.

2) There is only one element of our subset that is in $X_1\cup X_2\cup X_3\cup X_4$. We have got only two such subsets: $X_0\cup X_2$ and $X_0\cup X_3$.

All together we see that in our case there are $3n-1$ maximal cones.
Let us choose the following ordered sequence of $3n-1$ line bundles from $B_1\cup
B_2\cup B_3$:
\begin{equation}\label{list}
\begin{cases}
\o(-(n-1+b)D_v-D_y-D_t),\o(-(n-1+b)D_v-D_t),\\ \o(-(n-2+b)D_v-D_y-D_t),  
\o(-(n-2+b)D_v-D_t),... , \\
\o(-(b+1)D_v-D_t),\o(-bD_v-D_y-D_t),\\ \o(-(n-1)D_v),\o(-(n-2)D_v),...,\o.
\end{cases}
\end{equation}
We want to prove that this is a strongly exceptional collection.
We know that for any line bundles $L$ and $M$ on $X$:
$$Ext^i(L,M)=H^i(L^\vee\otimes M).$$ First we want to prove that for
any $L$ and $M$ in \eqref{list} $L^\vee\otimes M$ is acyclic. Let us write down line bundles of the form $L^\vee\otimes M$ 
where $L$ and $M$ are taken from \eqref{list}.
$$Diff=\begin{cases} (1)\;\;\;\o(sD_v) & s=-(n-1),...,n-1 \\
(2)\;\;\;\o(sD_v+D_t) & s=b+2-n,...,n-1+b \\
(3)\;\;\;\o(sD_v-D_t) & s=-(n-1+b),...,n-b-2 \\
(4)\;\;\;\o(sD_v+D_y) & s=-(n-1),...,n-2 \\
(5)\;\;\;\o(sD_v-D_y) & s=-(n-2),...,n-1 \\
(6)\;\;\;\o(sD_v+D_y+D_t) & s=b-(n-1),...,b+n-1 \\
(7)\;\;\;\o(sD_v-D_y-D_t) & s=-(b+n-1),...,n-1-b \\
\end{cases}$$
From Corollary \ref{acyclic} we know that they are acyclic if they are not
of the form
$$\o(\alpha_1D_{v}+\alpha_2D_{y}+\alpha_3D_{z}+\alpha_4D_{t}+\alpha_5D_{u})\cong$$
$$\cong \o((\alpha_1+\alpha_5+\alpha_3b)D_{v}+(\alpha_2+\alpha_3-\alpha_5)D_{y}+(\alpha_3+\alpha_4)D_{t})$$
where exactly $2,3$ or $5$ consecutive $\alpha$ are negative and
if $\alpha_1<0$ then $\alpha_1\leq-(n-1)$.

We will show that all line bundles of $Diff$ are not of this form. First
let us observe that they are not of this form for all $\alpha$
negative since then the coefficient of $D_t$ is less than or equal
to $-2$. Let us suppose that they are of this form with exactly
$2$ or $3$ consecutive $\alpha$ negative.

\bigskip (1) The coefficient of $D_y$ is $0$ therefore
$\alpha_2+\alpha_3=\alpha_5 $. But $\alpha_2,\alpha_3$ and
$\alpha_5$ cannot have the same sign (we treat $0$ as positive) so
$\alpha_2$ and $\alpha_3$ have different signs.  This means that
$\alpha_3$ and $\alpha_4$ have the same sign, and as
$\alpha_3+\alpha_4=0$, they both have to be equal to zero. So
$\alpha_2$ and as a consequence $\alpha_1$ are negative hence the
coefficient of $D_v$ is less than or equal to $-n$, which is a contradiction.

\bigskip (2) The coefficient of $D_y$ is $0$
so as before $\alpha_2$ and $\alpha_3$ have different signs. This
means that $\alpha_3$ and $\alpha_4$ are of the same sign. We know
that $\alpha_3+\alpha_4=1$, so they both have to be positive and
at most one equal to one. So $\alpha_2$ and as a consequence
$\alpha_1$ is negative hence the coefficient of $D_v$ is less than
or equal to $-(n-1)+b$, which is a contradiction.

\bigskip (3) As before $\alpha_2$ and $\alpha_3$ have different signs.
$\alpha_3$ cannot be positive since then $\alpha_4$ and as a consequence
coefficient of $D_t$ would also be positive. So $\alpha_3$ and as
a consequence $\alpha_4$ is negative hence the coefficient of
$D_t$ is less than or equal to $-2$, which is a contradiction.

\bigskip (4) The coefficient of $D_y$ is $1$ therefore
$\alpha_2+\alpha_3=\alpha_5+1$. But $\alpha_2,\alpha_3$ and
$\alpha_5$ cannot have the same sign,
so $\alpha_2$ and $\alpha_3$ have different signs or
$\alpha_2=\alpha_3=0$ and $\alpha_5=-1$.

First case: $\alpha_2$ and $\alpha_3$ have different signs. This means that
$\alpha_3$ and $\alpha_4$ have the same sign, and as $\alpha_3+\alpha_4=0$ we see that
$\alpha_3=\alpha_4=0$. So $\alpha_2$ and as a consequence
$\alpha_1$ is negative hence the coefficient of $D_v$ is less than
or equal to $-n-1$, which is a contradiction.

Second case: $\alpha_2=\alpha_3=0$ and $\alpha_5=-1$. We have also
$\alpha_4=0$ so $\alpha_1$ and $\alpha_5$ are negative hence the
coefficient of $D_v$ is less than or equal to $-n$, which is a contradiction.

\bigskip (5) The coefficient of $D_y$ is $-1$ therefore
$\alpha_2+\alpha_3=\alpha_5-1$. But $\alpha_2,\alpha_3$ and
$\alpha_5$ cannot have the same sign
so $\alpha_2$ and $\alpha_3$ have different signs.
This means that
$\alpha_3$ and $\alpha_4$ have the same sign, and as $\alpha_3+\alpha_4=0$ we see that
$\alpha_3=\alpha_4=0$.
 So $\alpha_2$ and as a consequence
$\alpha_1$ is negative hence the coefficient of $D_v$ is less than
or equal to $-(n-1)$, which is a contradiction.

\bigskip (6) The coefficient of $D_y$ is $1$ therefore
$\alpha_2+\alpha_3=\alpha_5+1$. But $\alpha_2,\alpha_3$ and
$\alpha_5$ cannot have the same sign,
so $\alpha_2$ and $\alpha_3$ have different signs or
$\alpha_2=\alpha_3=0$ and $\alpha_5=-1$.

First case: Assume that $\alpha_2$ and $\alpha_3$ have different signs. In this case $\alpha_3$  and $\alpha_4$ are of the same sign and as
$\alpha_3+\alpha_4=1$ they have to be positive and at most one. So
$\alpha_2$ and as a consequence $\alpha_1$ is negative hence the
coefficient of $D_v$ is less than or equal to $-n+b$, which is a
contradiction.

Second case: Assume $\alpha_2=\alpha_3=0$.
We have
$\alpha_4=1$ so $\alpha_1$ and $\alpha_5$ are negative hence the
coefficient of $D_v$ is less than or equal to $-n$, which is a contradiction.

\bigskip (7) The coefficient of $D_y$ is $-1$ therefore
$\alpha_2+\alpha_3=\alpha_5-1$. But $\alpha_2,\alpha_3$ and
$\alpha_5$ cannot have the same sign so $\alpha_2$ and $\alpha_3$
have different signs. $\alpha_3$ cannot be positive since then
$\alpha_4$ and as a consequence the coefficient of $D_t$ would also be
positive. So $\alpha_3$ and as a consequence $\alpha_4$ is
negative. Hence the coefficient of $D_t$ is less than or equal to
$-2$, which is a contradiction.

\bigskip
To have a strongly exceptional collection it remains to prove
that any pair of two line bundles $L_i$ and $L_j$ for $i<j$ from
our ordered sequence satisfies $0=Ext^0(L_j,L_i)=H^0(L_j^\vee\otimes
L_i)$. This is equivalent to showing that $L_j^\vee\otimes L_i$
has no global sections so from \cite{bohu} it is not of the form
$$(*)\hskip 6pt \o(\alpha_1D_{v}+\alpha_2D_{y}+\alpha_3D_{z}+\alpha_4D_{t}+\alpha_5D_{u})\cong$$
$$\cong \o((\alpha_1+\alpha_5+\alpha_3b)D_{v}+(\alpha_2+\alpha_3-\alpha_5)D_{y}+(\alpha_3+\alpha_4)D_{t})$$
with all $\alpha_i$ nonnegative. Let us partition our ordered collection into two collections:
$$Col_1=(\o(-(n-1+b)D_v-D_y-D_t),\o(-(n-1+b)D_v-D_t),$$$$\o(-(n-2+b)D_v-D_y-D_t),\o(-(n-2+b)D_v-D_t),$$$$...
,\o(-(b+1)D_v-D_t),\o(-bD_v-D_y-D_t))$$
and
$$Col_2=(\o(-(n-1)D_v),\o(-(n-2)D_v),...,\o).$$
If we take a difference of an element from $Col_1$ and $Col_2$, then the coefficient of $D_t$ is negative so the difference is not of the form (*). If we take the difference of two elements from $Col_2$ then the coefficient of $D_v$ is negative, hence it is not of the form (*). If we take the difference of two elements from $Col_1$, then either the coefficients of $D_v$ is negative or the difference is equal to $-D_y$. The divisor $-D_y$ is not of the form (*), because $\alpha_5$ would have to be strictly positive, hence the coefficient of $D_v$ would not be zero.

We have proven:

\ob\label{kol} With the above notation the following ordered sequence of line bundles in $X$
$$\o(-(n-1+b)D_v-D_y-D_t),\o(-(n-1+b)D_v-D_t),...$$
$$...,\o(-(b+1)D_v-D_t),\o(-bD_v-D_y-D_t),\o(-(n-1)D_v),...,\o$$
is a strongly exceptional collection.\kob \kwadrat

\subsection{Generating a derived category}

Finally, we will prove that the strongly exceptional collection given in Proposition \ref{kol} is also full. As already mentioned it is enough to
prove that it generates all line bundles of the set $B_1\cup B_2\cup
B_3$. In order to show that we need following two lemmas.

\lem\label{jeden} Let $k$ be any integer. Line bundles
$\mathcal{O}(-kD_v-D_y-D_t),\dots,$ $\mathcal{O}(-(n-1+k)D_v-D_y-D_t),\mathcal{O}(-(k+1)D_v-D_t),\dots, \mathcal{O}(-(n-1+k)D_v-D_t)$
generate $\mathcal{O}(-kD_v-D_t)$ in the derived category.\klem
\dow We consider the Koszul complex for
$\mathcal{O}(D_{y}),\mathcal{O}(D_{v_1}),\dots,\mathcal{O}(D_{v_{n-1}})$:
$$0\rightarrow \mathcal{O}(-(n-1)D_v-D_y)\rightarrow\dots\rightarrow\mathcal{O}(-D_v)^{n-1}\oplus\mathcal{O}(-D_y)\rightarrow\mathcal{O}\rightarrow 0.$$
By tensoring it with $\mathcal{O}(-kD_v-D_t)$ we obtain:
$$0\rightarrow \mathcal{O}(-(n-1+k)D_v-D_y-D_t)\rightarrow\dots\rightarrow$$
$$\rightarrow\mathcal{O}(-(1+k)D_v-D_t)^{n-1}\oplus\mathcal{O}(-k D_v-D_y-D_t)\rightarrow\mathcal{O}(-k D_v-D_t)\rightarrow 0.$$
All sheaves that appear in this exact sequence, apart from the
last one, are exactly
$\mathcal{O}(-kD_v-D_y-D_t),\dots,\mathcal{O}(-(n-1+k)D_v-D_y-D_t),\mathcal{O}(-(k+1)D_v-D_t),\dots,\mathcal{O}(-(n-1+k)D_v-D_t)$,
so indeed we can generate $\mathcal{O}(-k D_v-D_t)$. \kdow

\lem\label{dwa} Let $k$ be any integer. Line bundles
$\mathcal{O}(-(k+1)D_v-D_y-D_t),\dots,\mathcal{O}(-(n-1+k)D_v-D_y-D_t),\mathcal{O}(-(k+1)D_v-D_t),\dots,\mathcal{O}(-(n+k)D_v-D_t)$
generate $\mathcal{O}(-kD_v-D_y-D_t)$ in the derived
category.\klem \dow The proof is similar to the last one. We have
to consider the Koszul complex for line bundles
$\mathcal{O}(D_{u}),\mathcal{O}(D_{v_1}),\dots,\mathcal{O}(D_{v_{n-1}})$:
$$0\rightarrow \mathcal{O}(-(n-1)D_v-D_u)\rightarrow\dots\rightarrow\mathcal{O}(-D_v)^{n-1}\oplus\mathcal{O}(-D_u)\rightarrow\mathcal{O}\rightarrow 0$$
we dualize it and we tensor by $\mathcal{O}(-(n+k)D_v-D_t)$.\kdow
Summarizing, we have proved:
\tw\label{glowne} Let $X$ be a smooth, complete, $n$ dimensional toric variety with Picard number three, ray generators $X_0\cup\dots\cup X_4$, where
$$X_0=\{v_1,\dots,v_{n-1}\},\quad X_1=\{y\},\quad X_2=\{z\},\quad X_3=\{t\},\quad X_4=\{u\},$$
primitive collections $X_0\cup X_1$, $X_1\cup X_2,\dots,X_4\cup X_0$ and primitive relations:
$$v_1+\dots+v_{n-1}+y-(b+1)t=0,$$
$$y+z-u=0,$$
$$z+t=0,$$
$$t+u-y=0,$$
$$u+v_1+\dots+v_{n-1}-bt=0,$$
where $b$ is a positive integer.

Then the collection
$$\o(-(n-1+b)D_v-D_y-D_t),\o(-(n-1+b)D_v-D_t),...$$
$$...,\o(-(b+1)D_v-D_t),\o(-bD_v-D_y-D_t),\o(-(n-1)D_v),...,\o$$
is a full strongly exceptional collection of line bundles. \ktw
\dow We already know that this is a strongly exceptional
collection from Proposition \ref{kol}. Inductively using lemmas
\ref{jeden} and \ref{dwa} we can prove that it generates sets
$B_1$ and $B_2$. The set $B_3$ is already in our collection, hence
our collection is also full. \kdow

Arijit Dey

{Max-Planck Institut fur Mathematik, 

Vivatsgasse 7, 53111 Bonn, Deutchland}

e-mail address:\emph{arijit@mpim-bonn.mpg.de} 
\vskip 20pt

Micha\l \hskip 5pt Laso\'{n}

{Mathematical Institute of the Polish Academy of Sciences,  

\'{S}w. Tomasza 30, 31-027 Krak\'{o}w, Poland}
\vskip 2pt
{Theoretical Computer Science Department, 

Faculty of Mathematics and Computer Science, 

Jagiellonian University, 30-387 Krak\'{o}w, Poland}

e-mail address:\emph{mlason@op.pl}
\vskip 15pt

Mateusz Micha\l ek

{Mathematical Institute of the Polish Academy of Sciences,  

\'{S}w. Tomasza 30, 31-027 Krak\'{o}w, Poland}
\vskip 2pt
{Institut Fourier, Universite Joseph Fourier, 

100 rue des Maths, BP 74, 38402 St Martin d'He`res, France}

e-mail address:\emph{wajcha2@poczta.onet.pl}
\end{document}